\documentclass[11pt]{amsart}
\usepackage{amsmath,amsthm,amssymb}
\font\tenmsb=msbm10
\font\sevenmsb=msbm7
\font\fivemsb=msbm5

\newfam\msbfam
\textfont\msbfam=\tenmsb
\scriptfont\msbfam=\sevenmsb
\scriptscriptfont\msbfam=\fivemsb

\def\Bbb#1{\fam\msbfam\relax#1}

\newtheorem{thm}{Theorem}[section]
\newtheorem{prop}[thm]{Proposition}
\newtheorem{cor}[thm]{Corollary}
\newtheorem{lem}[thm]{Lemma}
\newtheorem{clm}{Claim}

\theoremstyle{definition}
\newtheorem{exa}[thm]{Example}
\newtheorem{defn}[thm]{Definition}
\newtheorem{note}[thm]{Notation}
\newtheorem{hyp}[thm]{Standing Hypothesis}

\theoremstyle{remark}
\newtheorem{rem}[thm]{Remark}
\newtheorem{comm}[thm]{Comment}
\newtheorem{openq}[thm]{Open Question}

\newcommand{\ben}{\begin{enumerate}}
\newcommand{\een}{\end{enumerate}}
\newcommand{\blem}{\begin{lem}}
\newcommand{\elem}{\end{lem}}
\newcommand{\bcl}{\begin{clm}}
\newcommand{\ecl}{\end{clm}}
\newcommand{\bthm}{\begin{thm}}
\newcommand{\ethm}{\end{thm}}
\newcommand{\bpr}{\begin{prop}}
\newcommand{\epr}{\end{prop}}
\newcommand{\bco}{\begin{cor}}
\newcommand{\eco}{\end{cor}}
\newcommand{\bcon}{\begin{conj}}
\newcommand{\econ}{\end{conj}}
\newcommand{\bde}{\begin{defn}}
\newcommand{\ede}{\end{defn}}
\newcommand{\bex}{\begin{exa}}
\newcommand{\eexa}{\end{exa}}
\newcommand{\bexe}{\begin{exe}}
\newcommand{\eexe}{\end{exe}}
\newcommand{\bobs}{\begin{obs}}
\newcommand{\eobs}{\end{obs}}
\newcommand{\bnote}{\begin{note}}
\newcommand{\enote}{\end{note}}
\newcommand{\beq}{\begin{equation}}
\newcommand{\eeq}{\end{equation}}
\newcommand{\bcom}{\begin{comm}}
\newcommand{\ecom}{\end{comm}}
\newcommand{\bhyp}{\begin{hyp}}
\newcommand{\ehyp}{\end{hyp}}
\newcommand{\brem}{\begin{rem}}
\newcommand{\erem}{\end{rem}}
\newcommand{\bopenq}{\begin{openq}}
\newcommand{\eopenq}{\end{openq}}
\numberwithin{equation}{section}

\newcommand{\Z}{{\Bbb Z}}
\newcommand{\C}{{\Bbb C}}
\newcommand{\Real}{{\Bbb R}}
\newcommand{\Q}{{\Bbb Q}}

\newcommand{\N}{{\Bbb N}}

\newcommand{\A}{\forall}

\newcommand{\z}{\zeta}
\newcommand{\x}{\xi}

\renewcommand{\a}{\alpha}
\newcommand{\g}{\gamma}
\newcommand{\be}{\beta}

\renewcommand{\L}{\Lambda}
\newcommand{\m}{\mu}

\renewcommand{\th}{\theta}

\newcommand{\R}{{\mathcal{R}}}
\renewcommand{\H}{{\mathcal{H}}}

\newcommand{\abs}[1]{\left\vert#1\right\vert}
\newcommand{\set}[1]{\left\{#1\right\}}
\newcommand{\sset}[1]{\bigl\{#1\bigr\}}
\newcommand{\bra}[1]{\left ( #1 \right )}
\newcommand{\Bra}[1]{\left [ #1 \right ]}
\def\comment#1{}

\begin{document}

\title{Generalized Widder Theorem via fractional moments}

\thanks{This paper is a part of the author's Ph.D. thesis,
written under the direction of Prof. Shmuel Kantorovitz in the
Department of Mathematics, Bar-Ilan University,
Israel.}

\author{Ami Viselter}
\address{Department of Mathematics, Bar Ilan University, 52900 Ramat-Gan, Israel}
\curraddr{Department of Mathematical and Statistical Sciences, University of Alberta, Edmonton, Alberta T6G 2G1, Canada}
\email{viselter@ualberta.ca}

\begin{abstract}
We provide a necessary and sufficient condition for the
representability of a function as the classical multidimensional
Laplace transform, when the support of the representing measure is
contained in some generalized semi-algebraic set. This is done by
employing a method of Putinar and Vasilescu [Putinar, M. and
Vasilescu, F.-H., Solving moment problems by dimensional extension,
Ann. of Math. (2) 149 (1999), no. 3, 1087-1107] for the
corresponding multidimensional moment problem.
\end{abstract}
\maketitle

\section*{Introduction}

A well-known theorem of Widder states that \emph{a necessary and
sufficient condition for a function $f : \Real \to \Real$ to be
representable in the form $$(\A x \in \Real) \quad f(x) =
\int_{-\infty}^\infty e^{-xt}d\m(t),$$ where $\m$ is a positive
measure over $\Real$, is that $f(x)$ be continuous and of positive
type} (c.f. \cite[Ch. VI, \S 21]{W}). This theorem has been
generalized to the multidimensional case in several works, using
various methods (see Akhiezer \cite{Akh}, Devinatz \cite{Dev},
Shucker \cite{Shucker} and their references, to mention a few), and
has many applications. If, however, we wished to characterize those
functions with a representing measure whose \emph{support} is
contained in some rather general, fixed set, $S$, the method used in
these papers would not be helpful, except for the case when $S$ is a
multidimensional box.

On the other hand, Putinar and Vasilescu used in \cite{PV} a method
of dimensional extension to solve the multidimensional moment
problem. In their paper, the moment problem is translated to the
problem of representation of a certain linear functional, which is
obtained by means of the spectral theory of selfadjoint operators,
over some special Hilbert space. The bonus in their method is that
it enables them to characterize moment sequences, whose representing
measure's support lies in a given semi-algebraic set.

The connection between a moment problem and the corresponding
Laplace transform representation problem has been successfully
established in the past (e.g. in \cite{W_doubly}, and the references
therein). In this note we modify Putinar and Vasilescu's method of
dimensional extension to obtain a generalized version of Widder's
Theorem, which characterizes the functions that can be represented
by the multidimensional Laplace transform of a measure with support
in a given (generalized) semi-algebraic set. Essentially,
non-negative integral powers of the variables are replaced by
non-negative rational powers.

\section{Preliminaries} \label{section:preliminaries}

Let $\R$ be an algebra of complex functions, such that $\overline{f}
\in \R$ for all $f \in \R$ (that is, $\R$ is selfadjoint). We say
that a linear functional $\L$ over $\R$ is \emph{positive
semi-definite} if $\L(\abs{f}^2) \geq 0$ for each $f \in \R$. When
this is the case, one can define the semi-inner product $(f,g) :=
\L(f\overline{g})$. Thus, if $\mathcal{N} = \sset{f \in \R :
\L(\abs{f}^2) = 0},$ then $\R / \mathcal{N}$ is an inner-product
space. Hence, its completion, $\H$, is a complex Hilbert space. For
simplicity, we often write $r$ instead of $r + \mathcal{N}$ for
elements $r \in \R / \mathcal{N}$.

The standard notations $\Real_+ = [0,\infty)$, $\Q_+ = \Real_+ \cap
\Q$, etc. are used. Fix an $n \in \N$. For $t=(t_1,\ldots,t_n), \a =
(\a_1,\ldots,\a_n) \in \Real_+^n$, we write $t^\a$ for
$t_1^{\a_1}\cdots t_n^{\a_n}$.

We let $\mathcal{P}_n$ denote the set of all complex polynomials
with $n$ real variables. By $\mathcal{Q}_n$ we shall denote the
complex algebra of all "fractional polynomials" of positive rational
exponents and $n$ variables. That is, $\mathcal{Q}_n$ is the set of
all of the functions in the form $\Real_+^n \ni t \mapsto \sum_{\a
\in \Q_+^n} a_\a t^\a$, where the $a_\a$'s are complex, and differ
from zero only for a finite number of indices $\a$.

Let $\Bbb A$ be a subsemigroup of $\Q_+^n$. A family of complex
numbers $\delta = (\delta_\a)_{\Bbb A}$ induces the linear
functional $L_\delta$ over the subalgebra of $\mathcal{Q}_n$
generated by $\set{t^\a : \a \in \Bbb A}$, defined by
$L_\delta(t^\a) = \delta_\a$ for all $\a \in \Bbb A$. We say that
$\delta$ is positive semi-definite if the functional $L_\delta$ is
positive semi-definite.

For the rest of the section, $\H$ denotes an arbitrary complex
Hilbert space.

\blem \label{frac_oper} Let $A$ be a positive selfadjoint operator
over $\H$, and let $q_1, q_2$ be positive real numbers. Then there
exists a \emph{unique} positive selfadjoint operator $B$, so that
$B^{q_2} = A^{q_1}$, namely $B = A^{q_1/q_2}$. \elem

\proof Let $E(\cdot; A)$ denote the resolution of the identity
of the selfadjoint operator A. By \cite[Theorem XII.2.9]{DS}, a
positive operator $B$ satisfies the theorem's statement if and only
if for every Borel set $\delta \subseteq \Real_+,$ \beq
\label{frac_oper_eq1} E(\delta^{1/q_2}; B) = E(\delta; B^{q_2}) =
E(\delta; A^{q_1}) = E(\delta^{1/q_1}; A). \eeq Since the mapping
$\delta \mapsto \delta^{q_2}$ is bijective from the set of all Borel
subsets of $\Real_+$ into itself, (\ref{frac_oper_eq1}) is
equivalent to that \beq \label{frac_oper_eq2} E(\delta; B) =
E(\delta^{q_2 / q_1}; A) \eeq for all Borel sets $\delta \subseteq
\Real_+$. But by the same theorem from \cite{DS}, there exists a
unique positive selfadjoint operator $B$ that satisfies
(\ref{frac_oper_eq2}), which is $B = A^{q_1 / q_2}$. \qed

We now state two results from \cite{PV}.

\bpr[{\cite[Proposition 2.1]{PV}}, originally from \cite{Nel, Nus}]
\label{VP_pr_2.1} Let $T_1, \ldots, T_n$ be symmetric operators in
$\H$. Assume that there exist a dense linear space $\mathcal{D}
\subseteq \cap_{j,k=1}^n D(T_jT_k)$ such that $T_jT_kx = T_kT_jx$
for all $x \in \mathcal{D}$, $j \neq k$, $j,k=1,\ldots,n$. If the
operator $(T_1^2+\cdots+T_n^2)_{|\mathcal{D}}$ is essentially
selfadjoint, then the operators $T_1,\ldots,T_n$ are essentially
selfadjoint, and their canonical closures $\overline{T_1},\ldots,
\overline{T_n}$ commute. \epr

\blem[{\cite[Lemma 2.2]{PV}}] \label{VP_lem_2.2} Let $A$ be a
positive densely defined operator in $\H$, such that $AD(A)
\subseteq D(A)$. Suppose that $I+A$ is bijective on $D(A)$. Then $A$
is essentially selfadjoint. \elem

\section{Generalized Widder Theorem} \label{section:main}

Let $p=(p_1,\ldots,p_m)$, where $p_k$ are real fractional
polynomials in $\mathcal{Q}_n$. For this fixed set of polynomials,
let $\th_p : \Real_+^n \to \C$ be defined as
$$\th_p(t) :=
(1+t_1^2+\ldots+t_n^2+p_1(t)^2+\ldots+p_m(t)^2)^{-1}.$$ We denote by
$\R$ the complex algebra generated by $\mathcal{Q}_n$ and the
function $\th_p$.

The following is the main operator-theoretic result, leading to the
moments theorem to follow.

\bthm \label{frac_mom_op_thm} Let $\L$ be a positive semi-definite
functional over $\R$. Then there exists a unique representing
measure for $\L$. The support of that measure is contained in
$\Real_+^n$. Moreover, if $\L(p_k \abs{r}^2) \geq 0$ for all $r \in
\R$, $1 \leq k \leq m$, then the support of that (unique) measure is
a subset of $\bigcap_{k=1}^m p_k^{-1}(\Real_+)$. \ethm

\proof Let $\H$ be the Hilbert generated by $\L$, as explained in
Section \ref{section:preliminaries}. For $1 \leq i \leq n$, $1 \leq
j \leq m$, we define the operators $T_i, P_j$ over $\R /
\mathcal{N}$ by $$T_i : r + \mathcal{N}\mapsto t_i r + \mathcal{N},
\quad P_j : r + \mathcal{N} \mapsto p_j r + \mathcal{N}.$$ Let $B$
be the operator $B := T_1^2 + \ldots T_n^2 + P_1^2 + \ldots +
P_m^2$. Then  $B : \R / \mathcal{N} \to \R / \mathcal{N}$ is a
positive operator, since for all $r \in \R$, $(Br,r) = \sum_{i=1}^n
\L(\abs{t_i r}^2) + \sum_{j=1}^m \L(\abs{p_j r}^2) \geq 0$, by the
positivity of $\L$. Moreover, $I+B$ is bijective, since for all $r
\in \R$, $(I+B)u = r$ for some $u \in \R$ if and only if $u = \th_p
r$. Therefore, by Lemma \ref{VP_lem_2.2}, $B$ is essentially
selfadjoint. Thus, by Proposition \ref{VP_pr_2.1}, the operators
$T_i$ and $P_j$ are essentially selfadjoint for all $1 \leq i \leq
n$, $1 \leq j \leq m$. Moreover, the selfadjoint operators $A_1 :=
\overline{T_1},\ldots,A_n := \overline{T_n}$ commute, and thus have
a common resolution of the identity, $E$ (c.f. \cite[Ch. IV, Theorem
10.3]{V}). Set $A := (A_1,\ldots,A_n)$. For $r \in \R$, $r(A)$ will
denote the normal operator $\int_{\Real^n_+} r(t)E(dt)$, and for $q
\in \Q_+^n$, $A^q$ will stand for $f(A)$ where $f(t) =t^q$. We
define the operators $T_i(q_i) : \R / \mathcal{N} \to \R /
\mathcal{N}$, $1 \leq i \leq n$, by
$$T_i(q_i) : r + \mathcal{N} \mapsto t_i^{q_i} r + \mathcal{N},$$ and set $T(q) :=
T_1(q_1)\cdots T_n(q_n)$.

\bcl \label{claim_1} For all $q \in \Q_+^n$, $T(q) \subseteq A^q$.
\ecl To prove the claim, we first notice that $T(q)$ is positive for
$q \in \Q_+^n$, as $\L(t^q \abs{r}^2) = \L(\abs{t^{q/2}r}^2) \geq 0$
for all $r \in \R$. Let
$q=(\frac{k_1}{\ell_1},\ldots,\frac{k_n}{\ell_n})$, where
$k_1,\ldots,k_n \in \N \cup \set{0} ,\ell_1,\ldots,\ell_n \in \N$.
Fix an $1 \leq i \leq n$. Since the operator $T_i(\frac{1}{\ell_i})$
is positive, it has an (a priori, not necessarily unique) positive
selfadjoint extension, $A_i(\frac{1}{\ell_i})$. Now, observe that
$T_i = T_i(\frac{1}{\ell_i})^{\ell_i} \subseteq
A_i(\frac{1}{\ell_i})^{\ell_i}.$ But $T_i$ is essentially
selfadjoint and $A_i(\frac{1}{\ell_i})^{\ell_i}$ is selfadjoint,
which implies that $A_i = \overline{T_i} =
A_i(\frac{1}{\ell_i})^{\ell_i}$. Therefore, by the uniqueness part
of Lemma \ref{frac_oper}, $ T_i(\frac{1}{\ell_i}) \subseteq
A_i(\frac{1}{\ell_i}) = A_i^{1/\ell_i}$. Hence, once again by Lemma
\ref{frac_oper}, and the fact that
\begin{equation}\label{eq:mult_op_oper_calc} (\A r_1, r_2 \in \R)
\quad r_1(A)r_2(A) \subseteq (r_1 r_2)(A)\end{equation} (which
follows readily from \cite[Ch. IV, Theorem 10.3]{V}),
\begin{equation}\begin{split}T(q) & = T_1(\frac{1}{\ell_1})^{k_1} \cdots
T_n(\frac{1}{\ell_n})^{k_n} \subseteq \\ & \subseteq (A_1^{1/\ell_1})^{k_1} \cdots
(A_n^{1/\ell_n})^{k_n} = A_1^{k_1 / \ell_1} \cdots
A_n^{k_n / \ell_n} \subseteq A^q,\end{split}\end{equation} and the claim is proved.

\bcl For all $r \in \R$, \beq \label{rep} \L(r) =
\int\limits_{\Real_+^n} r(t) \bra{E(dt)(1+\mathcal{N}),
1+\mathcal{N}}.\eeq \ecl In order to prove the claim, fix an $r \in
\R$. Let the operator $r(T) : \R / \mathcal{N} \to \R / \mathcal{N}$
be the operator of multiplication by $r$. We shall show that $r(T)
\subseteq r(A)$. By linearity and \eqref{eq:mult_op_oper_calc}, it
is sufficient to prove this for $r(t) = t^q$, $q \in \Q_+^n$, and
for $r=\th_p$. The first case is exactly Claim 1, since $r(T) =
T(q)$ and $r(A) = A^q$. As for the case $r=\th_p$, it follows from
the fact that $\th_p^{-1}(T) \subseteq \th_p^{-1}(A)$, and so for
all $f \in \R / \mathcal{N}$, $\th_p(A)f =
\th_p(A)\Bra{\th_p^{-1}(T)\th_p(T)}f =
\th_p(A)\th_p^{-1}(A)\th_p(T)f = \th_p(T)f$ (by
\eqref{eq:mult_op_oper_calc}). Finally, to prove \eqref{rep}, we
note that by the Spectral Theorem,
$$\L(r) = (r+\mathcal{N}, 1+\mathcal{N}) = \bra{r(T)(1+\mathcal{N}),
1+\mathcal{N}} = \bra{r(A)(1+\mathcal{N}), 1+\mathcal{N}} = $$
$$ = \int\limits_{\Real_+^n} r(t) \bra{E(dt)(1+\mathcal{N}),
1+\mathcal{N}}$$ (the domain of integration is $\Real_+^n$ since the
operators $A_1,\ldots,A_n$ are positive), and the claim is proved.

Consequentially, the (positive) Borel measure $\m$ over $\Real^n$
defined by $\m(\cdot) = (E(\cdot)(1+\mathcal{N}), 1+\mathcal{N})$ is
a representing measure for $\L$, whose support lies in $\Real_+^n$.
We have thus proved the existence part of the theorem.

The uniqueness of $\mu$ is proved as in \cite{PV}, using an argument
taken from \cite{Fug}. Let us assume that there exists another
positive measure, $\nu$, over $\Real^n_+$, that represents the
functional $\L$. It is clear that when this is the case, $\R /
\mathcal{N}$ can be identified as a subspace of the Hilbert space
$L^2(\nu)$. Hence, $\mathcal{H}$ can be identified as the closure of
$\R / \mathcal{N}$ in $L^2(\nu)$. For all $1 \leq j \leq n$, let us
now define the selfadjoint operators $H_j$ over $L^2(\nu)$ by $H_j f
:= t_j f.$ Denote the spectral measure of $H_j$ by $E_j$. Since the
operators $H_1, \ldots, H_n$ commute, they have a joint spectral
measure, $E_H$. Obviously, $T_j \subseteq H_j$ for all $j$. Since
the operators $H_j$ are closed, $A_j \subseteq H_j$ for all $j$.
Therefore, $R(\z; A_j) \subseteq R(\z; H_j)$ for all $\z \in \C
\backslash \Real$, and so $R(\z; H_j)$ leaves $\mathcal{H}$
invariant, whence we conclude (c.f. \cite[Theorem XII.2.10]{DS})
that $E_j$ also leaves $\mathcal{H}$ invariant for each $1 \leq j
\leq n$. Thus, as $E_H(B_1 \times \ldots \times B_n) = E_1(B_1)
\cdots E_n(B_n)$ for all Borel sets $B_1, \ldots, B_n$ in $\Real$,
$E_H$ leaves $\mathcal{H}$ invariant as well. In particular, for
each Borel set $B$ in $\Real^n$, $I_B = E_H(B)1 \in \mathcal{H}$
(where $I_B$ is the indicator function of $B$ over $\Real^n$). Since
the simple functions are dense in $L^2(\nu)$, we infer that
$\mathcal{H} = L^2(\nu)$, and so $A_j = H_j$ for each $1 \leq j \leq
n$. In particular, $E = E_H$, and so for each Borel set $B$ in
$\Real^n$, by the definition of $\m$,
$$\m(B) = (E(B)(1+\mathcal{N}), 1+\mathcal{N}) = (E_H(B)1, 1) =
\int\limits_{\Real^n} I_B d\nu = \nu(B),$$ and the proof of the
uniqueness of $\m$ is completed.

Assume now that $\L(p_k \abs{r}^2) \geq 0$ for all $r \in \R$ and $1
\leq k \leq m$. This condition is equivalent to the operators $P_1,
\ldots, P_m$ being positive. We recall that these operators are
essentially selfadjoint. But for all such $k$, $P_k \subseteq
p_k(A)$ by Claim \ref{claim_1}, and $p_k(A)$ is selfadjoint; thus,
$\overline{P_k} = p_k(A)$ is a positive selfadjoint operator.
Equivalently, its spectral measure is supported by $\Real_+$. But
the spectral measure of $p_k(A)$ is $F_k(\delta) =
E(p_k^{-1}(\delta))$. Hence, $E$ itself is supported by
$p_k^{-1}(\Real_+)$. Since that is true for all $1 \leq k \leq m$,
the support of $E$ is therefore a subset of $\bigcap_{k=1}^m
p_k^{-1}(\Real_+)$. \qed

\blem \label{vas_put_lemma_2.3} Let $\vartheta \in \mathcal{P}_n$ be
such that $\vartheta(t)
> 0$ for all $t \in \Real^n$, and let $p(t,s) \in \mathcal{P}_{n+1}$
($t \in \Real^n$, $s \in \Real$) be such that $p(t,
\vartheta^{-1}(t)) \equiv 0$. Then there exists a complex polynomial
$q \in \mathcal{P}_{n+1}$ such that $$(\A t,s) \quad p(t,s) = q(t,s)
\cdot \Bra{s\vartheta(t)-1}.$$ \elem

\proof This is a simple generalization of \cite[Lemma 2.3]{PV};
simply replace their $\theta_\textbf{p}$ by $\vartheta^{-1}$. We
omit the details.

\bde Denote by $\widetilde{\mathcal{Q}}_n$ the complex algebra
generated by $\mathcal{Q}_n$ and the algebra of all complex
polynomials with one positive real variable. Its elements will take
the form $p(t,s)$, $t \in \Real^n_+$, $s \in \Real_+$.\ede

\bpr \label{alg_iso} Let $\rho : \widetilde{\mathcal{Q}}_n \to \R$
be the mapping defined by $p(t,s) \mapsto p(t, \th_p(t))$. Then
$\rho$ is a surjective algebras homomorphism, whose kernel is the
ideal generated by the function
$$\sigma(t,s) = s\th_p(t)^{-1}-1.$$ \epr

\proof We first note that $\rho$ is indeed well-defined, since
$\th_p(t) \in \Real_+$ for all $t \in \Real^n_+$. It is clearly a
surjective algebras homomorphism. Assume $p \in \ker(\rho)$, that
is, $p(t, \th_p(t)) = 0$ for all $t \in \Real^n_+$. For $1 \leq j
\leq n$, we let $c_j$ denote the l.c.m of all of the denominators of
the exponents of $t_j$ in the polynomial $p$. The mappings $u_j =
t_j^{1/c_j}$ are bijective mappings from $\Real_+$ onto itself.
Replacing $t_j$ by $u_j^{c_j}$ in the above equality yields \beq
\label{kernel_thm_eq} (\A u \in \Real^n_+) \quad p(u^c, \th_p(u^c))
= 0. \eeq The expression on the left side of (\ref{kernel_thm_eq}),
after the reduction of the fractions in the exponents of the
$u_j$'s, becomes a (not fractional) polynomial in $u = (u_1, \ldots,
u_n)$. Hence, (\ref{kernel_thm_eq}) is true (as equality of
polynomials) for all $u \in \Real^n$, and by Lemma
\ref{vas_put_lemma_2.3}, there exists a $q \in \mathcal{P}_{n+1}$
such that $$(\A u \in \Real^n_+, s \in \Real_+) \quad p(u^c, s) =
q(u,s) \Bra{s\th_p(u^c)^{-1}-1}.$$ We can now replace $u$ by
$t^{1/c}$, and by defining $\widetilde{q}(t,s) = q(t^{1/c}, s)$, we
conclude that
$$p(t,s) = \widetilde{q}(t,s) \Bra{s\th_p(t)^{-1}-1}$$ for all $t
\in \Real^n_+$, $s \in \Real_+$, where $\widetilde{q} \in
\widetilde{\mathcal{Q}}_n$, as wanted. \qed

\bde Let $\g = (\g_\a)_{\a \in \Real_+^n}$ be a family of
non-negative numbers. \ben
\item We say that $\g$ is continuous if the function $\a \mapsto
\g_\a$ is continuous (as a function from $\Real_+^n$ to $\Real_+$).
\item We say that $\g$ is an
\emph{($n$-dimensional) fractional moments family} if there exists a
positive Borel measure, $\m$, over $\Real_+^n$, such that
\begin{equation} \label{def:frac_mom_fam} (\A \a \in \Real_+^n) \quad \g_\a = \int_{\Real_+^n} t^\a
d\m.\end{equation}
 \een \ede

Note that \eqref{def:frac_mom_fam} is equivalent to the
(multidimensional) Laplace representation $$(\A \a \in \Real_+^n)
\quad \g_\a = \int_{\Real^n} e^{-\a \cdot s} d\nu(s)$$ obtained by
the change of variable $t = e^{-s}$.

The following is the main theorem, whose proof is almost identical
to that of Theorem 2.7 in \cite{PV}, basing on our Theorem
\ref{frac_mom_op_thm} and Proposition \ref{alg_iso} instead of the
parallel ones in \cite{PV}, and using Lebesgue's Dominated
Convergence Theorem to derive \eqref{def:frac_mom_fam} for all of
$\Real_+^n$. For the sake of completeness, we include the details.

\bthm \label{frac_mom_thm} Let $\g = (\g_\a)_{\a \in \Real_+^n}$ be
a continuous family of non-negative numbers. Let $p_1, \ldots, p_m
\in \mathcal{Q}_n$, $p_k(t) = \sum_{\x \in I_k} a_{k\x}t^\x$ ($I_k
\subseteq \Q_+^n$ is finite) for $k=1, 2, \ldots, m$. Then $\g$ is a
fractional moments family with a representing measure whose support
is a subset of $\cap_{k=1}^m p_k^{-1}(\Real_+)$ if and only if there
exists a positive semi-definite family $$\delta =
(\delta_{(\a,\be)})_{(\a,\be) \in \Q_+^n \times \Z_+}$$ that
satisfies: \ben
\item $\delta_{(\a,0)} = \g_\a$ for all $\a \in \Q_+^n$.
\item $\delta_{(\a,\be)} = \delta_{(\a, \be+1)} + \sum_{j=1}^n
\delta_{(\a+2e_j, \be+1)} + \sum_{k=1}^m \sum_{\x, \eta \in I_k}
a_{k\x} a_{k\eta} \delta_{(\a+\x+\eta, \be+1)}$ for all $(\a, \be)
\in \Q_+^n \times \Z_+$.
\item The families $\bra{\sum_{\x \in I_k}
a_{k\x}\delta_{(\a+\x,\be)}}_{(\a,\be) \in \Q_+^n \times \Z_+}$ are
positive semi-definite for all $k = 1,\ldots, m$. \een Moreover, the
representing measure of $\g$ (with the properties mentioned above)
is unique if and only if the family $\delta$ is unique. \ethm

\proof Necessity. Assume that $\g$ is a fractional moments family
with a representing measure, $\m$, whose support is a subset of $E
:= \cap_{k=1}^m p_k^{-1}(\Real_+)$. We define the family $\delta$ by
\beq \label{frac_mom_thm_eq1} (\A (\a,\be) \in \Q_+^n \times \Z_+)
\quad \delta_{(\a,\be)} := \int_E t^\a \th_p(t)^\be d\m. \eeq Then
$\delta$ is a positive semi-definite family, that satisfies (1). (2)
is a result of the obvious equality
$$\int_E \bra{\th_p(t)(1+t_1^2 + \ldots + t_n^2 + p_1(t)^2 + \ldots
p_m(t)^2) - 1} t^\a \th_p(t)^\be d\m = 0,$$ which is true for all
$\a \in \Q_+^n$, $\be \in \Z_+$. Finally, (3) is true since
$$\int_E p_k(t) \abs{p(t,\th_p(t))}^2 d\m \geq 0$$ for all $p \in
\widetilde{\mathcal{Q}}_n$, $1 \leq k \leq m$.

Sufficiency. Let $\delta$ be as in the theorem's statement, and the
algebra $\R$ be defined as in the beginning of this section. We
define the linear functional $\L$ over $\R$ by
$$\L(r) = L_\delta(p)$$ for all $r \in \R$, where $L_\delta$ is the
linear functional induced by $\delta$ over
$\widetilde{\mathcal{Q}}_n$, and $p \in \widetilde{\mathcal{Q}}_n$
is such that $r(t) = p(t, \th_p(t))$ for all $t \in \Real_+^n$. $\L$
is well-defined, since by Proposition \ref{alg_iso}, $\R \cong
\widetilde{\mathcal{Q}}_n / \mathcal{I}$, where $\mathcal{I}$ is the
ideal in $\widetilde{\mathcal{Q}}_n$, generated by the element
$s\th_p(t)^{-1}-1$; and indeed, by (2), $(L_\delta)_{| \mathcal{I}}
= 0$. Thus, $\L$ is a well-defined positive semi-definite mapping on
$\R$. From (3) we deduce that $L_\delta(p_k \abs{p}^2) \geq 0$ for
all $p \in \widetilde{\mathcal{Q}}_n$, $1 \leq k \leq m$, hence
$\L(p_k \abs{r}^2) \geq 0$ for all $r \in \R$, $1 \leq k \leq m$.

By virtue of Theorem \ref{frac_mom_op_thm}, there exists a unique
representing measure, $\m$, for $\L$, whose support is a subset of
$E$. Particularly, by (1), \beq \label{frac_mom_thm_eq2} \g_\a =
\delta_{(\a, 0)} = \int_E t^\a d\m \eeq for all $\a \in \Q_+^n$. But
since the family $\g$ is continuous, Lebesgue's Dominated
Convergence Theorem implies that $\g_\a = \int_E t^\a d\m$ for all
$\a \in \Real_+^n$, that is, $\g$ is a fractional moments sequence,
as wanted.

Assume that the family $\delta$, that satisfies the conditions in
the theorem's statement, is unique. Let $\m_1, \m_2$ be two
representing measures for $\g$. By the uniqueness of $\delta$,
equation (\ref{frac_mom_thm_eq1}) and the discussion that follows,
$\int_E t^\a \th_p(t)^\be d\m_1 $ $= \int_E t^\a \th_p(t)^\be d\m_2$
for each $\a \in \Q_+^n$, $\be \in \Z_+$. Therefore, $\int_E r d\m_1
= \int_E r d\m_2$ for all $r \in \R$, and by the uniqueness part of
Theorem \ref{frac_mom_op_thm}, it follows that $\m_1 = \m_2$.

Conversely, assume that $\m$ is unique. Suppose that both
$\delta_1$, $\delta_2$ satisfy the conditions in the theorem's
statement. As explained above, $\delta_1, \delta_2$ induce the
positive semi-definite linear functionals $\L_1, \L_2$,
respectively, over $\R$, which, in turn, have the unique
representing measures $\m_1, \m_2$, respectively (by Theorem
\ref{frac_mom_op_thm}). Both measures represent $\g$ as a fractional
moments family, and so, by the uniqueness of $\m$, $\m_1 = \m_2$,
hence $\L_1 = \L_2$. Finally, for each $\a \in \Q_+^n$, $\be \in
\Z_+$, $(\delta_1)_{(\a,\be)} = \L_1(t^\a \th_p(t)^\be) = \L_2(t^\a
\th_p(t)^\be) = (\delta_2)_{(\a,\be)}$, that is, $\delta_1 =
\delta_2$. \qed

\section*{Acknowledgments}
The author would like to thank his thesis advisor, Prof. Shmuel Kantorovitz, for his encouragement,
many valuable discussions, and for having suggested the topic of this paper.

\end{document}